\documentclass[11pt]{amsart}
\usepackage[left=1in,right=1in]{geometry}
\usepackage{amssymb,amsmath,amsthm}
\usepackage[flushmargin]{footmisc}
\usepackage[foot]{amsaddr}
\usepackage{wasysym}
\usepackage{graphicx}
\usepackage{epsfig,subfig}
\usepackage{color}
\usepackage[nocompress]{cite}

\title{Perturbation of a nonsmooth supercritical Hopf bifurcation}

\author{Julie Leifeld}

\begin{document}

\newtheorem{theorem}{Theorem}
\newcommand{\ode}[2]{\dfrac{d{#1}}{d{#2}}}
\newcommand{\pde}[2]{\frac{\partial{#1}}{\partial{#2}}}
\newcommand{\red}[1]{\textcolor{red}{#1}}
\newcommand{\sgn}{\mathop{\mathrm{sgn}}}
\newcommand{\Tr}{\text{Tr}}

\begin{abstract}

Nonsmooth formulations of physical models are common, particularly in climate modeling.  However, in many of these models, there is little justification for this modeling choice, and no mathematical indication that the resulting behavior in the nonsmooth system is consistent with  behavior that would arise in a similar smooth system.  I explore this question by analyzing the relationship between a smooth and nonsmooth Hopf bifurcation in one particular climate example.  I show that the nonsmooth bifurcation perturbs to the smooth bifurcation under small changes in a limiting parameter.

\end{abstract}

\maketitle

\section{Introduction}
\label{sec:intro}

Because oscillatory behavior is of interest in many climate applications, climate scientists often look for Hopf bifurcations in their models.  Moreover, many climate phenomena are well modeled by nonsmooth systems, moving the possible bifurcation analyses outside the realm of classical dynamics \cite{Stommel1961,Welander82,Dercole07,hill2015analysis,mcgehee2014quadratic}.  
Here we focus on one particular conceptual climate model, Welander's ocean convection model \cite{Welander82}.  This is an ocean circulation box model, which divides the ocean into two regions, the surface ocean and deep ocean.  Mixing between the boxes is governed by a piecewise defined function, to demonstrate an abrupt transition between mixing and nonmixing states.   This discontinuous function is a simplification of a corresponding smooth function which represents the ``real" ocean.  The model finds oscillatory behavior in both the smooth and nonsmooth versions, but one naturally wonders whether the bifurcation structure seen in the analysis of the nonsmooth system imitates the structure of the smooth counterpart.  This question has been discussed in a few different contexts.  Sotomayor-Teixeira regularization uses {\it{smoothing functions}} to create smooth systems which agree with the original nonsmooth system outside an epsilon region around a discontinuity \cite{teixeira,sotomayor1996regularization}.  However, even if the use of a particular smoothing function is proven to demonstrate the same qualitative behavior as the nonsmooth system, this does little to answer the question, if one is trying to ultimately draw conclusions about phenomena in a previously given smooth system.  Moreover, it has been shown that in general, standard nonsmooth analysis techniques might not give qualitatively similar phenomena to nearby smooth systems \cite{Jeffrey11}.  Instead, an exact knowledge of the behavior of the system on the {\it{splitting manifolds}}, or discontinuity boundaries, is necessary to gain a full understanding of a nonsmooth system.  If this behavior is not known, the only current alternative is to know the behavior in the original smooth system, which defeats the purpose of the nonsmooth simplification.  However, in certain models, such as Welander's model, the nonsmooth system can give some information about the behavior in the smooth system.  We discuss this relationship for a specific Hopf bifurcation structure in Welander's model.  In Section \ref{sec:analysisbkg} we discuss the terminology we use to analyze the nonsmooth version of the model.  In Section \ref{sec:model} we outline the convection model, and in Sections \ref{sec:analysis} and \ref{sec:blowup} we analyze the model in a nonsmooth setting, and prove the persistence of the bifurcation point as the system perturbs to a smooth system.  Examples of this kind may be a segment of the path toward a comprehensive theory about the relationship between smooth and nonsmooth systems.

\section{Analysis of Nonsmooth systems}
\label{sec:analysisbkg}

Before the model can be introduced and analyzed, we must clarify the nonsmooth syntax we will use, and highlight the uncertainties inherent in standard nonsmooth analysis techniques.  In general, a nonsmooth system, with a splitting manifold defined by the nullcline of a scalar function $h(x)$, can be written in the form
\begin{equation}
\label{eq:general}
\dot{x}=f\left(x;\lambda\right),
\end{equation}
where $f$ depends on the variable $x\in\mathbf{R}^n$ and a nonsmooth parameter
\[
\lambda=\left\{\begin{array}{cc}1&h(x)>0\\ 0 & h(x)<0.\end{array}\right.
\]
One can then define regions of the phase space in terms of $\lambda$, and say $f^+(x)=f(x;1)$, and $f^-(x)=f(x;0)$, with $f^+$ and $f^-$ smooth.

\vspace{1em}

\begin{figure}[t]
	\centering
	\includegraphics[width=3in]{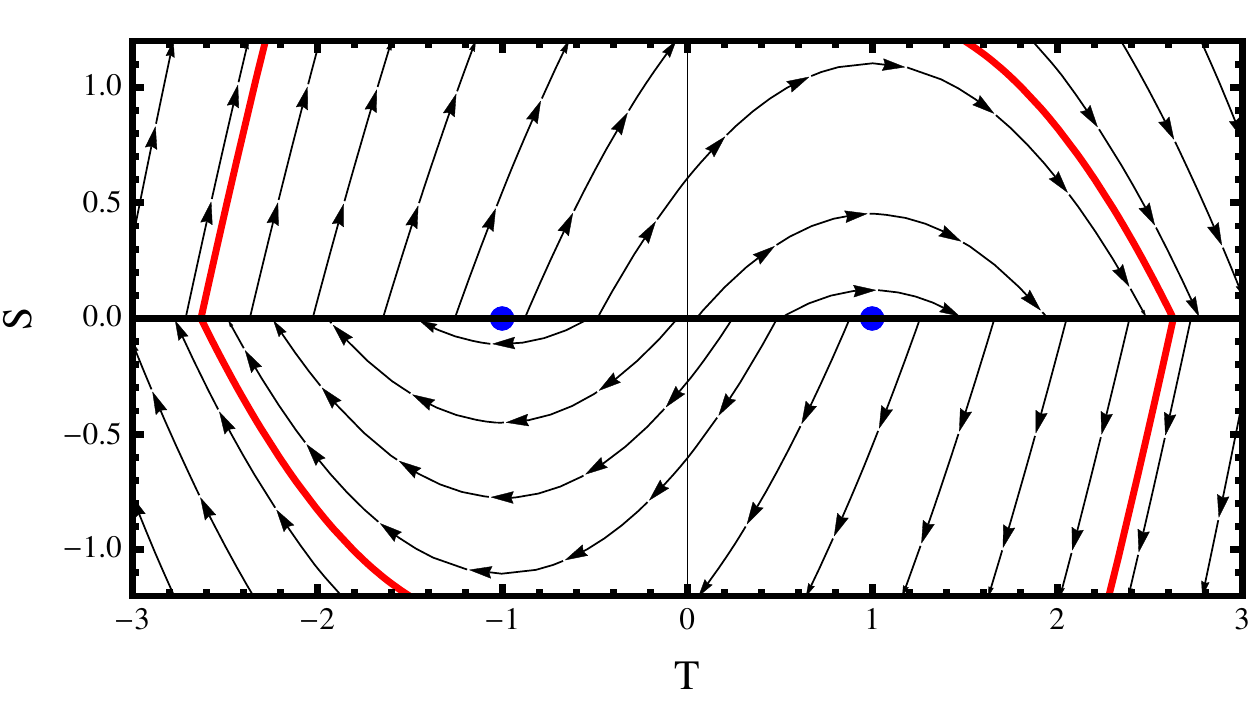}
	%\subfloat[]{\includegraphics[scale=1]{figures/ebmNonsmooth.eps}}
	%\subfloat[]{\includegraphics[width=0.45\textwidth]{figures/ebmFig.eps}}
	%\subfloat[]{\includegraphics[width=0.45\textwidth]{figures/ebmNonsmooth.eps}}
	\caption{A schematic example of an unstable sliding region.  Crossing solutions, like those in red, occur when the vector fields point in the same direction.  The boundaries of the sliding region occur at tangencies of the vector field to the splitting manifold, marked here by blue dots.}
	\label{fig:stableslide}
\end{figure}

In each smooth region of the phase space the system behaves in the manner dictated by classical smooth dynamical systems.  However, behavior on the splitting manifold is more complicated.  Solutions can either cross the manifold or slide along it, as in Figure \ref{fig:stableslide}.  Intuitively, one expects crossing if the vector fields on either side of the splitting manifold point in the same direction.  Alternatively, sliding regions occur when both vector fields point toward the splitting manifold.  Formally, sliding solutions exist in regions where the following system of equations can be solved for $\lambda\in[0,1]$:
\begin{equation}
\label{eq:s}
\begin{array}{r}S=f(x;\lambda)\cdot\nabla h(x)=0\\
h(x)=0. \end{array}
\end{equation}
Regions in which sliding solutions exist are unsurprisingly called {\it{sliding regions}}, and regions in which solutions must cross the splitting manifold are called {\it{crossing regions}}.  The stability of the sliding regions is given by the $\lambda$ derivative of $S$, with
\[
\ode{}{\lambda} S<0
\]
implying a stable sliding region, and 
\[
\ode{}{\lambda} S>0
\]
implying unstable sliding.

The flow in a sliding region can be found through the solution to \eqref{eq:s}.  For $\lambda^*$ solving \eqref{eq:s}, the flow in the sliding region is given by the system
\begin{equation}
\label{eq:slideflow}
\dot{x}=f(x;\lambda^*).
\end{equation}

\vspace{1em}

Solutions to \eqref{eq:general} are in general not unique, because many solutions can hit the sliding region in finite time, and merge.  In fact, one should not even assume that there is a unique solution $\lambda^*$ defining the sliding flow.  The number of solutions is a result of the dependence of $f$ on $\lambda$.  The standard and most commonly used method of analysis, formulated by Filippov \cite{filippov}, does not have this difficulty, as it is assumed that $f(x;\lambda)=\lambda f^+(x)+(1-\lambda)f^-(x)$.

\vspace{1em}

However, work by Jeffrey has shown that one might in some sense expect nonlinear dependence of $f$ on $\lambda$ in the form of $f(x;\lambda)=\lambda f^+(x)+(1-\lambda)f^-(x)+g(x;\lambda(1-\lambda))$ \cite{Jeffrey14hidden}.  These terms are problematic because they are only defined on the splitting manifold, and hence can drastically change behavior there, while leaving the behavior in the smooth regions unchanged.  This means that infinitely many smooth systems, each with qualitatively different behavior, can limit pointwise to the same nonsmooth system.  Moreover, nonsmooth systems that are the pointwise limits of functions of sigmoids, like Welander's model, have asymptotic expansions that contain these terms.  So, each nonsmooth system has infinitely many possible qualitative behaviors, and without previous knowledge of how the system behaves on the splitting manifold, the dependence of solutions on $f^+(x)$ and $f^-(x)$ is not clear.  

\vspace{1em}

There are several practical reasons for ignoring these difficulties, the most obvious being that assuming a linear dependence on $\lambda$ simplifies the analysis of the system.  Additionally, the real dependence of $f$ on $\lambda$ is not always clear in individual applications, and so in some cases this is the best approximation available.  However, this makes it clear that it is very important to understand the relationship between Filippov systems and the smooth systems that are perturbations of them, or alternatively the relationship between smooth systems and their limiting nonsmooth systems.  If a Filippov system is not able to capture all possible qualitative behaviors in the smooth systems that are nearby, then one should naturally ask which structures that are seen in nonsmooth systems are analogous to structures in the corresponding smooth system.  One side of this question is being addressed through the regularization method of Sotomayor and Teixera \cite{sotomayor1996regularization}.  In Sotomayor-Teixeira regularization, the splitting manifold in a Filippov system is smoothed in a strip of radius $\epsilon$, using a smoothing function.  Depending on one's choice of smoothing function, different qualitative behaviors can be seen in this region, as to be expected from the previous discussion.  Work has been done to give conditions on the smoothing function which result in the same structures in the smoothed system as in the Filippov system, for various nonsmooth phenomena \cite{buzzi2006singular,llibre2009study, kristiansen2015regularizations,teixeira}.  However, one must then also initially believe that the Filippov system shows the correct structures, which requires prior knowledge of both the smooth and the nonsmooth analysis in each specific application, and which is at the very least twice the work.  However, it is possible in some cases that one can determine the existence of structures in smooth systems simply by analyzing the Filippov limit of the system.  This paper uses a straightforward technique to show that the Hopf bifurcation in the Filippov analysis of the nonsmooth version of an ocean convection model is in fact also present in the corresponding smooth version of the model.

\section{Welander's Ocean Convection Model}
\label{sec:model}

\begin{figure}[t]
	\centering
	\includegraphics[height=3in]{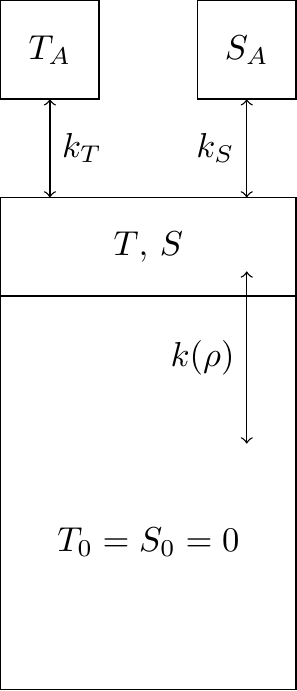}
	%\subfloat[]{\includegraphics[scale=1]{figures/ebmNonsmooth.eps}}
	%\subfloat[]{\includegraphics[width=0.45\textwidth]{figures/ebmFig.eps}}
	%\subfloat[]{\includegraphics[width=0.45\textwidth]{figures/ebmNonsmooth.eps}}
	\caption{A schematic diagram of Welander's box model.}
	\label{fig:welschematic}
\end{figure}

Welander's model \cite{Welander82} is a standard ocean box model, meaning that the ocean is divided into well mixed boxes, and circulation is controlled by density differences between the boxes.  The relevance of the model is the ability to find oscillatory behavior in a large scale, simple model, which both justified the use of simple models, and demonstrated the ocean's ability to display internally driven periodicity.  Welander's model separates the ocean into two boxes, a surface ocean and a deep ocean, as in Figure \ref{fig:welschematic}.  The relative density of the water in the boxes is controlled by the temperature and salinity of the water, through an equation of state: $\rho=-\alpha T+\beta S$.  The deep ocean is assumed to have constant density ($T_0=S_0=\rho_0=0$), as any changes in the temperature and salinity happen on a time scale that did not concern Welander.  Equations for the temperature and salinity in the surface box can then be written as 
\begin{equation}
\label{eq:Wel}
\begin{array}{ccl}
\dot{T} & = & k_T(T_A-T)-k(\rho)T\\[8pt]
\dot{S}&=&k_S(S_A-S)-k(\rho)S\\[8pt]
\rho&=&-\alpha T+\gamma S,
\end{array}
\end{equation}
where $\rho$ is the density of the surface box.  $T_A$ and $S_A$ are generic atmospheric temperature and salinity forcing, due to precipitation, evaporation, solar forcing, etc., and are assumed to be constant.  $k_T$ and $k_S$ are relaxation rates for the atmospheric forcing, and $k(\rho)$ is the rate of convective mixing, dependent on density.  We immediately nondimensionalize the model for ease of analysis, and will refer only to the nondimensionalized system for the rest of the paper.
\begin{equation}
\label{eq:Welandersystem}
\begin{array}{c}\begin{array}{rcl} \dot{T}&=&1-T-k(\rho)T\\ \dot{S}&=&\beta(1-S)-k(\rho)S, \end{array} \\[8pt]
 \\
\rho=-\alpha T+S.
\end{array}
\end{equation}

The key feature of Welander's model is the convective mixing function $k(\rho)$.  Welander makes several assumptions about this function, the most important being that convection is small when the density of the surface ocean is small, large when the density is large, and the transition between these states is abrupt.  These are reasonable assumptions, because the ocean is in general highly stratified, without much mixing between the layers.  It is the switching behavior in $k$ that allows for periodicity in Welander's model.  It is also this switching behavior that suggests that the system might be well modeled as nonsmooth.  In fact, Welander does a basic analysis for two different versions of the function $k$, one where $k$ is an inverse tangent, 
\begin{equation}
\label{eq:ksmooth}
k(\rho)=\frac{1}{\pi}\tan^{-1}\left(\frac{\rho-\varepsilon}{a}\right)+\frac{1}{2},
\end{equation}
and one where $k$ is a Heaviside function,
\begin{equation}
\label{eq:knonsmooth}
k(\rho)=\left\{\begin{array}{rl}1 & \rho>\varepsilon \\ 0 &\rho<\varepsilon, \end{array}\right.
\end{equation}
The nonsmooth version of the model is the limit of the smooth version, as the smoothness parameter $a\rightarrow0$.  In his paper, Welander does numerical simulations to show the existence of a periodic orbit in both the smooth and nonsmooth versions of his model.

\begin{figure}[t]
	\centering
	\includegraphics[width=3in]{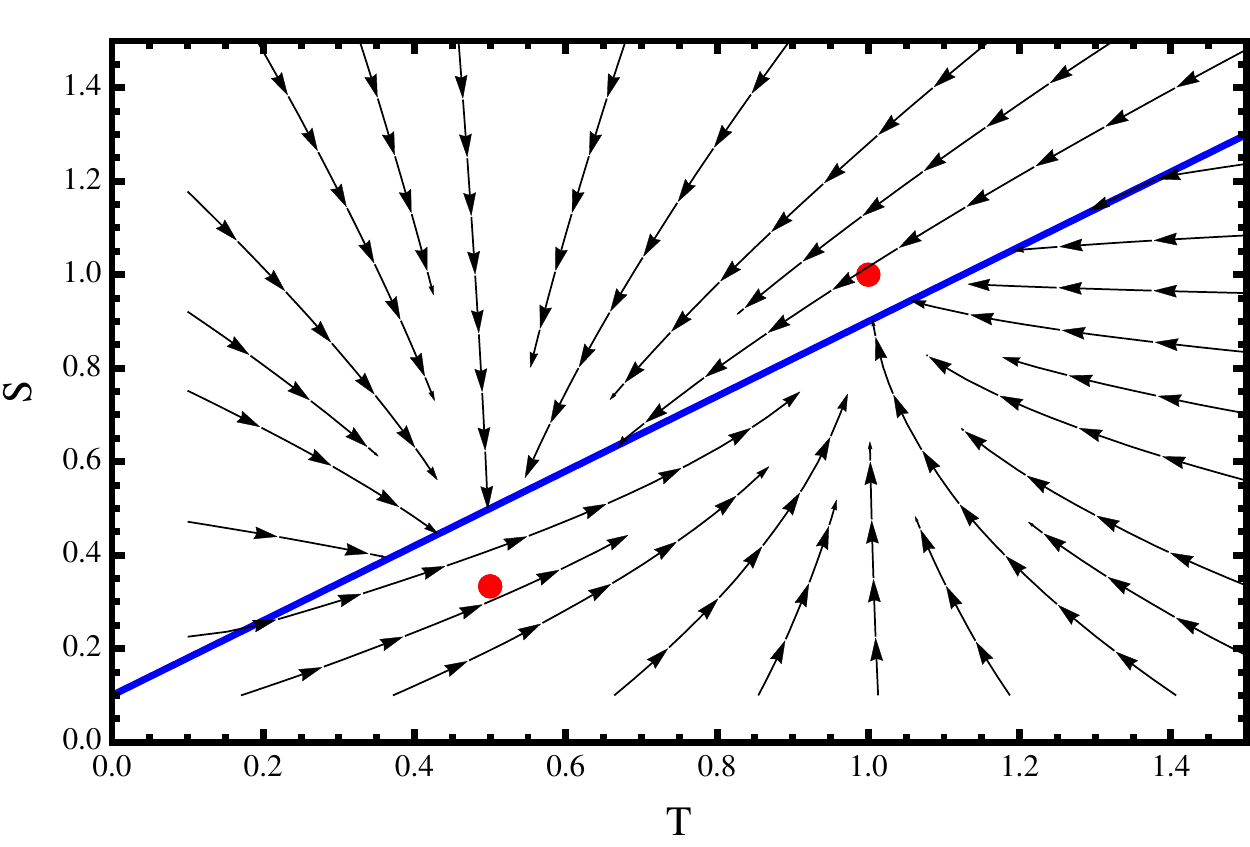}
	%\subfloat[]{\includegraphics[scale=1]{figures/ebmNonsmooth.eps}}
	%\subfloat[]{\includegraphics[width=0.45\textwidth]{figures/ebmFig.eps}}
	%\subfloat[]{\includegraphics[width=0.45\textwidth]{figures/ebmNonsmooth.eps}}
	\caption{A phase diagram for equation \eqref{eq:Welandersystem}, where $k$ is nonsmooth.  Virtual equilibria are given by red dots, and the splitting manifold is the blue line.}
	\label{fig:virteq}
\end{figure}

Welander chooses values for $\alpha$ and $\beta$, 
\begin{equation}
\label{alpha}
\alpha=\dfrac{4}{5},
\end{equation}
and
\begin{equation}
\label{beta}
\beta=\dfrac{1}{2}.
\end{equation}
Then, when Welander chooses $\varepsilon=-0.1$, oscillations in the nonsmooth model are due to two stable {\it{virtual equilibria}}, which means neither of the equilibria found by setting $k=0$ and $k=1$ in \eqref{eq:Welandersystem} exist in the region of phase space for which they satisfy their respective equations (see Figure \ref{fig:virteq}).  This leads to solutions that cross the splitting manifold, and then change directions to approach the other virtual equilibrium.  However, this type of oscillatory behavior does not guarantee the stable periodic orbit found by Welander.  In fact, for different values of $\varepsilon$, the stable periodic orbit does not exist, although the virtual equilibria remain virtual.  A standard analysis of the nonsmooth system indicates that the stable periodic orbit arises from a Hopf-like bifurcation, at a double tangency point which is called a {\it{fused focus}} \cite{filippov}.  This will be discussed in more detail in Section \ref{sec:analysis}.

\section{Filippov Analysis of Welander's Model}
\label{sec:analysis}

We start with a Filippov analysis of the Hopf bifurcation in the nondimensionalized nonsmooth version of Welander's model, \eqref{eq:Welandersystem},
\[
\begin{array}{rcl} \dot{T}&=&1-T-k(\rho)T\\ \dot{S}&=&\beta(1-S)-k(\rho)S, \end{array}
\]
\[
k(\rho)=\left\{\begin{array}{rl}1 & \rho>\varepsilon \\ 0 &\rho<\varepsilon, \end{array}\right.
\]
\[
\rho=-\alpha T+S.
\]
This is in the standard form seen in \eqref{eq:general}, where $\lambda=k$.  Because the splitting manifold, which is the line $\rho=-\alpha T+S=\varepsilon$, depends explicitly on the bifurcation parameter $\varepsilon$, it is useful to the analysis to do an initial change of coordinates,
\[
\begin{array}{rcl} x&=&T\\
y&=&S-\alpha T-\varepsilon, \end{array}
\]
after which the splitting manifold is the line $y=0$.  With some basic calculations, the new system is
\begin{equation}
\label{eq:Welcc}
\begin{array}{l}
\begin{array}{rcl}\dot{x}&=&1-x-k(y)x\\
\dot{y}&=&\beta-\beta\varepsilon-k\varepsilon-\alpha-(\beta+k)y-(\alpha\beta-\alpha)x
\end{array}\\
 \\
k(y)=\left\{\begin{array}{rl}1 & y>0\\ 0 & y<0\end{array}\right. .
\end{array}
\end{equation}

It is important to note that this coordinate change is a linear transformation composed with a translation.  Neither the translation or the linear transformation change the Filippov analysis of the nonsmooth system.  This is due to the fact that Filippov analysis depends on solving equations that are linear in the parameter $\lambda$, or $k$ in this system.

Now, the boundaries of the sliding region are given by points where $\dot{y}=0$ along $y=0$, which gives the equation
\[
0=\beta+\beta\varepsilon-k\varepsilon-\alpha-(\alpha\beta-\alpha)x.
\]
Solving this equation for $x$, and plugging in the values for $\alpha$ and $\beta$ given in \eqref{alpha} and \eqref{beta}, and $k=0$ or $k=1$, we find that the sliding region has boundaries 
\[
x=\frac{3}{4}+\frac{15}{4}\varepsilon,
\]
and 
\[
x=\frac{3}{4}+\frac{5}{4}\varepsilon.
\]

The first thing to notice is that the sliding region collapses to a single point in the case $\varepsilon=0$.  This point is called a {\it{double tangency}}.

\begin{figure}[t]
	\centering
	\subfloat[]{\includegraphics[width=0.45\textwidth]{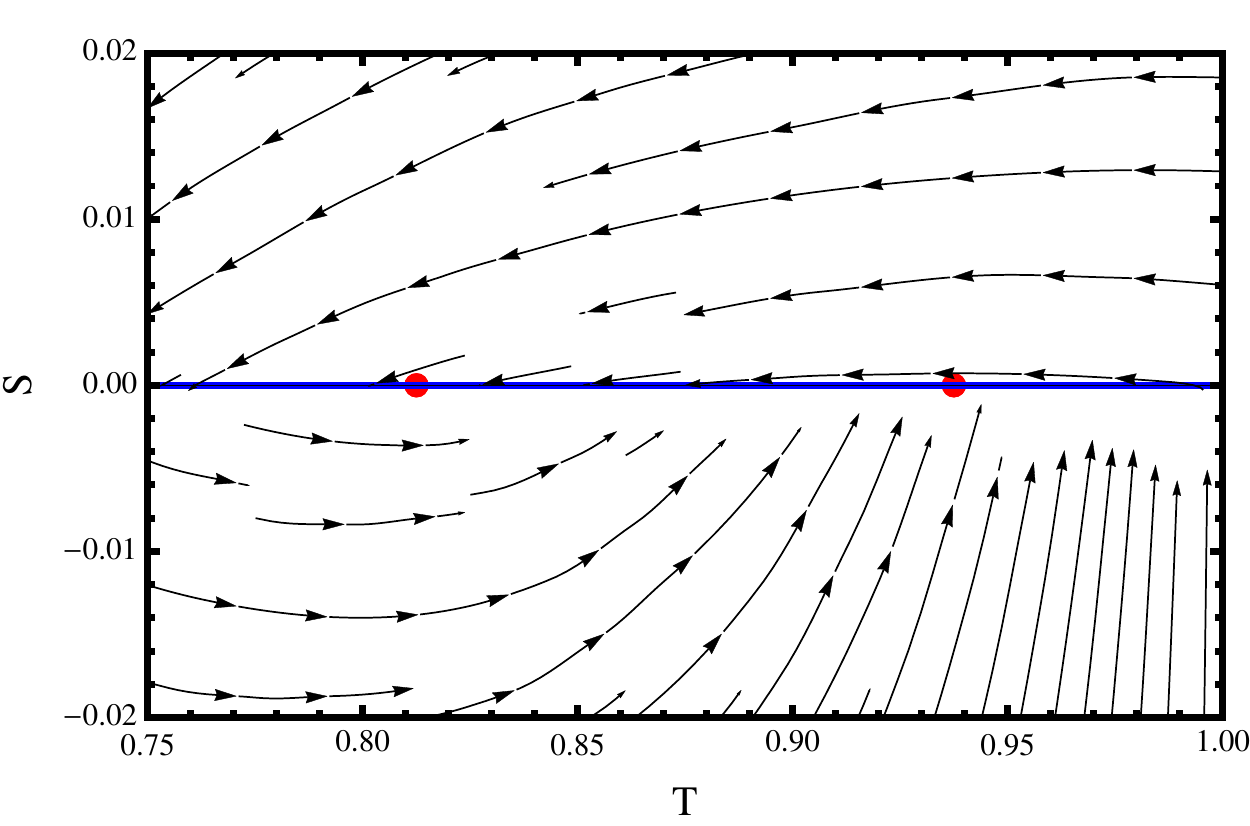}}
	\subfloat[]{\includegraphics[width=0.45\textwidth]{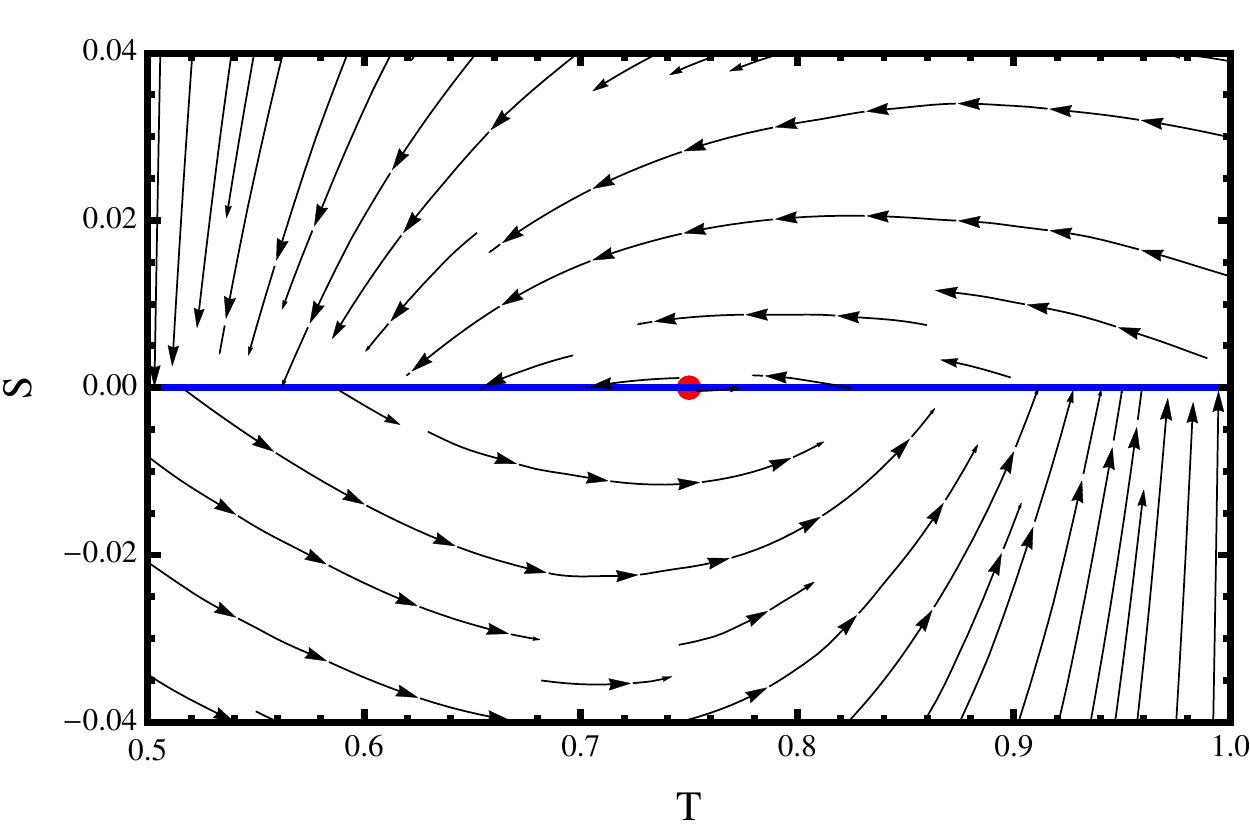}}
	
	\subfloat[]{\includegraphics[width=0.45\textwidth]{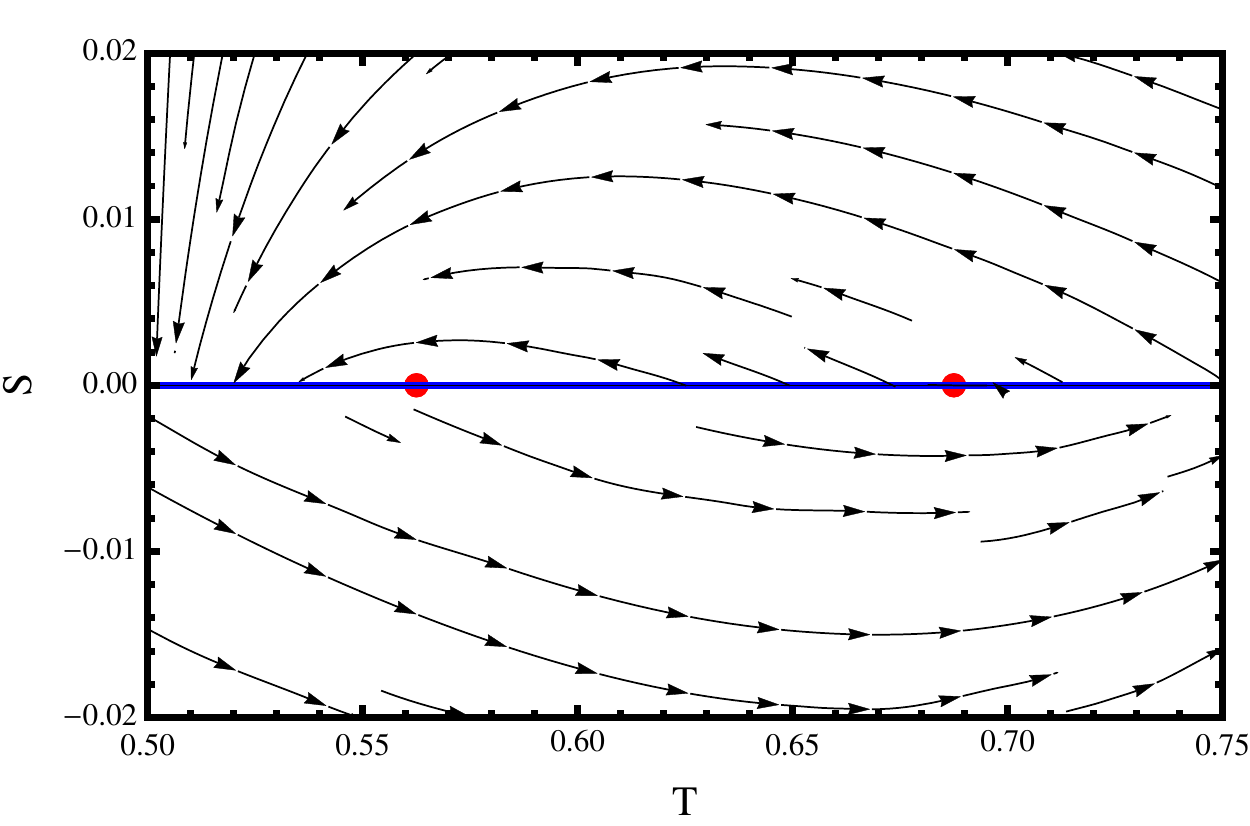}}
	\caption{The stability of the sliding region in system \eqref{eq:Welcc}.  In (A) $\varepsilon=0.05$ and the sliding region is stable. In (B) there is a double tangency when $\varepsilon=0$, and in (C) $\varepsilon=-0.05$, and the sliding region is unstable.}
	\label{fig:slidingstability}
\end{figure}

Stability of the sliding region is given by \eqref{eq:s}.  In the context of Welander's model, $\lambda=k$, and because of the coordinate change, we can choose $h=y$, which will satisfy all of the stipulations laid out in Section \ref{sec:analysisbkg}.  Here, 
\[
\ode{}{\lambda} (f\cdot\nabla h)|_{y=0}=\ode{}{k}\left(\beta-\beta\varepsilon-k\varepsilon-\alpha-(\beta+k)y-(\alpha\beta-\alpha)x\right)|_{y=0}=-\varepsilon.
\]
This indicates that the sliding region is stable for $\varepsilon>0$ and unstable for $\varepsilon<0$ (see Figure \ref{fig:slidingstability}).  It is reasonable that this change in stability occurs simultaneously to the sliding region narrowing to a point before expanding again.  This transition is also analogous to the change in stability necessary in a smooth Hopf bifurcation.  However, here there is no real equilibrium.   Instead there is a double invisible tangency, which Filippov called a {\it{fused focus}} \cite{filippov}.

The fused focus bifurcation is dependent on both tangencies being {\it{invisible}} when $\varepsilon=0$, meaning that the trajectory through the tangency does not exist in that region of phase space.  This allows for the double tangency to be a focus, as opposed to a saddle or some mixture of the two.  It is easy to see that the double tangency is invisible.  The $x$ derivative of the phase curves along $y=0$ is
\[
\ode{y}{x}=\dfrac{\beta-\beta\varepsilon-k\varepsilon-\alpha-(\alpha\beta-\alpha)x}{1-x-kx}.
\]
In the region where $k=1$, this function has a local maximum at the tangency, and when $k=0$, it has a local minimum.

\begin{figure}[t]
	\centering
	\subfloat[]{\includegraphics[width=0.45\textwidth]{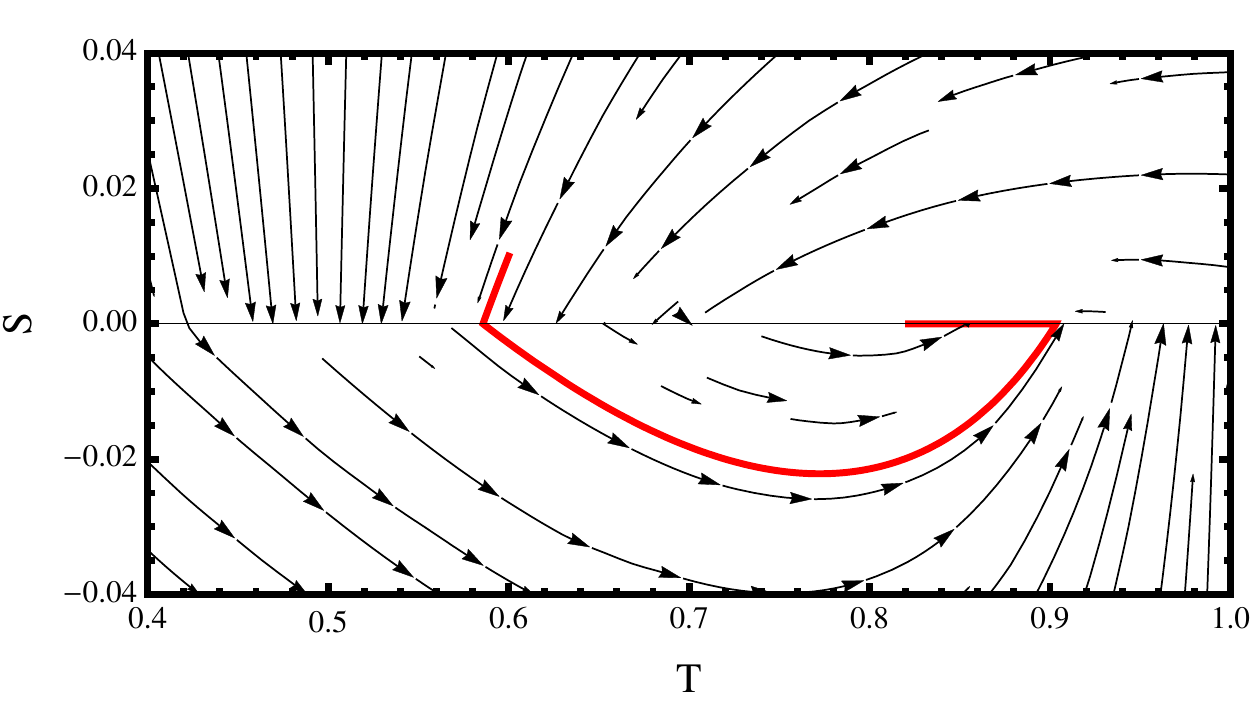}}
	\subfloat[]{\includegraphics[width=0.45\textwidth]{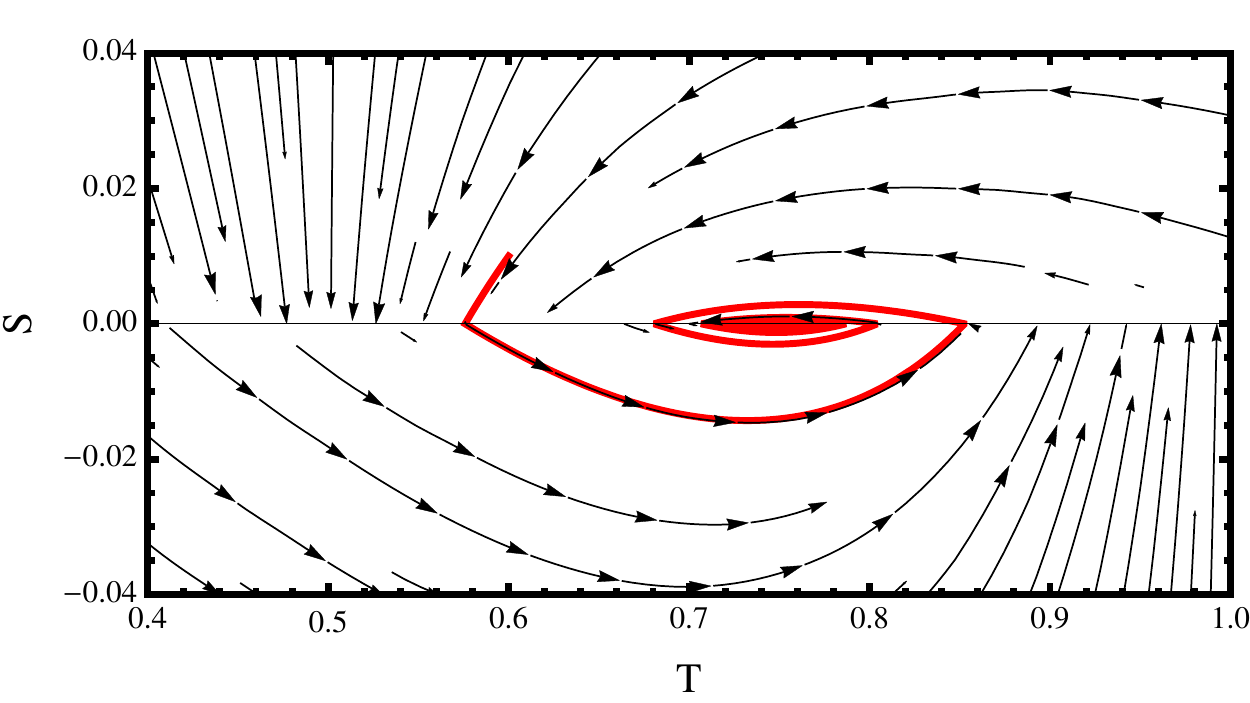}}
	
	\subfloat[]{\includegraphics[width=0.45\textwidth]{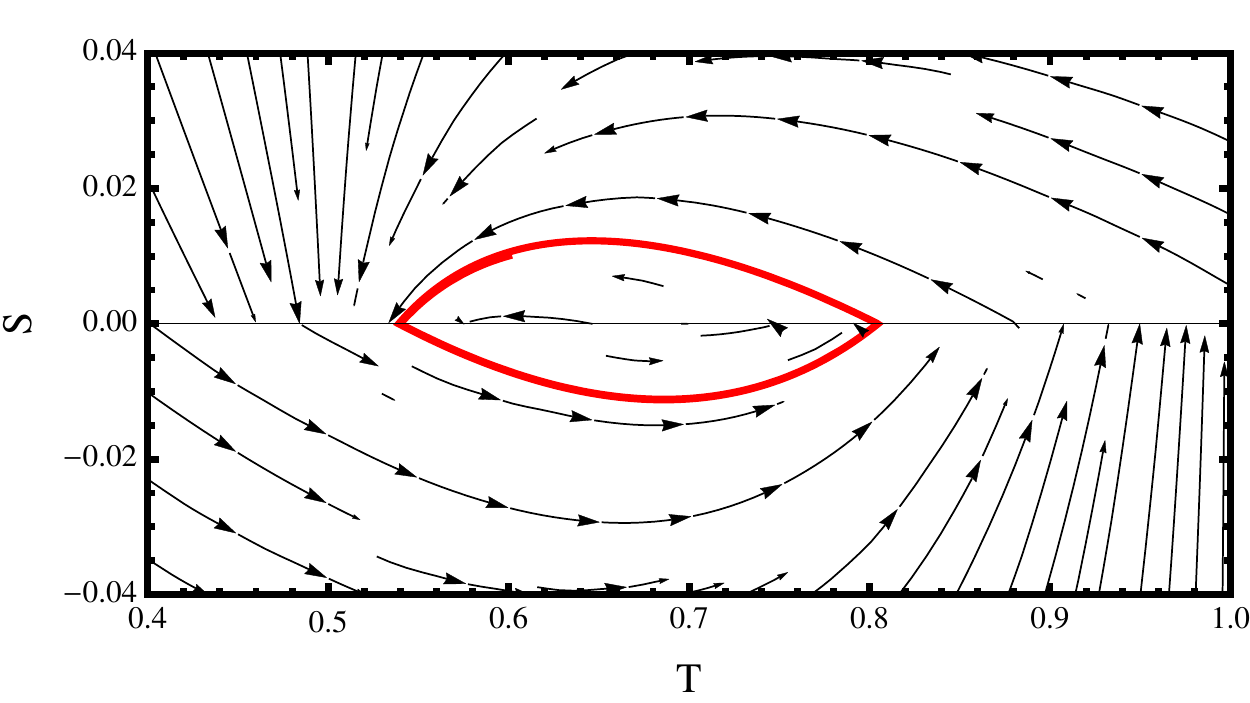}}
	\caption{The birth of a periodic orbit through the fused focus bifurcation.  (A) When $\varepsilon=0.04$, trajectories converge to the stable sliding region.  (B) When $\varepsilon=0$, trajectories converge to a stable focus.  (C) When $\varepsilon=-0.04$, trajectories converge to a stable periodic orbit.}
	\label{fig:periodicorbit}
\end{figure}

Filippov gave conditions on the stability of the fused focus \cite{filippov}.  In another long calculation, it can be shown that the fused focus at $\varepsilon=0$ is stable.  Then, the changing stability of the sliding region as $\varepsilon$ becomes negative immediately gives the existence of a small periodic orbit nearby.  The uniqueness of this orbit is not immediately clear, although simulations indicate that it is likely to be unique.  Figure \ref{fig:periodicorbit} shows trajectories of the system as it undergoes the fused focus bifurcation.

So, we see a nonsmooth version of a supercritical Hopf bifurcation using a Filippov analysis of Welander's model.

\section{Blow Up Method}
\label{sec:blowup}

Because the fused focus bifurcation in Welander's system occurs along the splitting manifold $y=0$, in the regular Filippov system it is impossible to do any standard analysis.  So, it is useful to do a coordinate change to instead look at the system in terms of independent variables $x$ and $k$.  We call this the ``blow up system," because it is effectively zooming in on the behavior of the system near $y=0$.  It should be explicitly noted that this is different from the blow up methods of Sotomayor-Teixera regularization, in which the system is smoothed in a size $\delta$ region around the splitting manifold.  Instead, this coordinate changed system is equivalent to the original system, as long as the smoothness parameter $a\neq0$.  The method of using a coordinate change to transform a singularity into a regular point has been used in several other mathematical areas, notably celestial mechanics \cite{mcgehee1974triple}.

Here, using \eqref{eq:ksmooth}, if 
\[
z=\dfrac{y}{a},
\]
we see that 
\begin{equation}
\label{eq:kwithz}
k=\Phi(z)=\frac{1}{\pi}\tan^{-1}(z)+\frac{1}{2}.
\end{equation}
The chain rule tells us
\[
\begin{array}{rcl}\dot{k}&=\ode{\Phi}{z}\ode{z}{y}\ode{y}{t}&=\dfrac{1}{a}\Phi'\left(\Phi^{-1}(k)\right)\left(\beta-\beta\varepsilon-k\varepsilon-\alpha-a(\beta+k)\Phi^{-1}(k)-(\alpha\beta-\alpha)x\right).
\end{array}
\]
For the $k$ as defined, 
\[
\Phi^{-1}(k)=\tan\left(\pi(k-\frac{1}{2})\right)=-\cot\left(\pi k\right)
\]
and
\[
\Phi'(y)=\frac{1}{\pi}\frac{1}{1+z^2}.
\]
Therefore the blow up system in terms of $x$ and $k$ is
\begin{equation}
\label{eq:Welblowup}
\begin{array}{rcl}
\dot{x}&=&1-x-kx\\
\dot{k}&=& \dfrac{1}{a\pi}\sin^2(\pi k)\left(\beta-\beta\varepsilon-k\varepsilon-\alpha+a(\beta+k)\left(\cot(\pi k)\right)-(\alpha\beta-\alpha)x\right)
\end{array}
\end{equation}
There are a few important things to note about this system.  First, for all $a\neq0$ it is equivalent to the original system.  The places where difficulties arise are where the function $\Phi$ is not invertible, at $k=0$ and $k=1$.  The derivative $\Phi'(z)$ is bounded everywhere, as the dependence on the parameter $a$ has been removed by the coordinate change.  Therefore, this system is well defined and smooth everywhere away from the $k$-boundaries.  

An equilibrium of the Filippov flow exists under the condition that the system of equations
\[
\begin{array}{c}
\beta-\beta\varepsilon-k\varepsilon-\alpha-(\alpha\beta-\alpha)x=0\\
1-x-kx=0
\end{array}
\]
is solvable for a $k\in[0,1]$.  For $\varepsilon$ near $\varepsilon=0$, the bifurcation point, this system is solvable, although for $|\varepsilon|$ sufficiently large, this pseudoequilibrium disappears.

Showing that the Filippov equilibrium perturbs to the smooth equilibrium is a singular perturbation problem if one uses the blow up system \eqref{eq:Welblowup}.  For an introduction to geometric singular perturbation theory, or GSPT, the reader is referred to the work of Jones \cite{Jones1995}.

Because $a$ is a small parameter, the blow up system is fast slow, with $k$ being the fast variable and $x$ being the slow variable, and one can analyze this system using GSPT techniques.  The system gives a critical manifold of 
\[
\beta-\beta\varepsilon-k\varepsilon-\alpha-(\alpha\beta-\alpha)x=0.
\]  
An equilibrium is a point on this critical manifold for which 
\[
x=\frac{1}{1+k}.  
\]
These are also exactly the conditions for a Filippov equilibrium in the nonsmooth $(x,y)$ system.  As long as $\varepsilon\neq0$, this critical manifold is immediately seen to be normally hyperbolic, so singular perturbation theory gives the existence for small $a>0$ of a unique equilibrium in the neighborhood of the Filippov equilibrium.  An application of the Intermediate Value Theorem shows that the equilibrium in the smooth version of Welander's model is unique.  Therefore, the Filippov equilibrium perturbs to the unique smooth equilibrium.  The stability of the equilibrium changes when $\varepsilon$ goes through 0.

However, GSPT breaks down when $\varepsilon=0$, because the critical manifold is no longer normally hyperbolic.  In fact, the blow up system gives a line of equilibria, along 
\[
x=\frac{3}{4}.
\]  
This is expected, as the bifurcation point is also only center stable.  However, the blow up system gives a unique candidate for a point around which to analyze the system.  This is the intersection of the vertical line 
\[
x=\frac{3}{4}
\]
and the nullcline 
\[
x=\frac{1}{1+k},
\]
which gives the relevant point
\begin{equation}
\label{eq:relpoint}
(x,k,\varepsilon,a)=\left(\frac{3}{4},\frac{1}{3},0,0\right).  
\end{equation}
Because the system is smooth away from $k=0$ and $k=1$, the equilibrium will limit to this point as $a\rightarrow0$ and $\varepsilon\rightarrow0$.  It should be noted that this is a two dimensional limit, so generally having independent limits as these parameters approach zero is not good enough.  However, the smoothness of the system allows this conclusion.

There is a Hopf bifurcation point in a neighborhood of the fused focus.  The fast system, 
\begin{equation}
\label{eq:blowupfast}
\begin{array}{rcl}
x'&=&a(1-x-kx)\\
k'&=& \dfrac{1}{\pi}\sin^2(\pi k)\left(\beta-\beta\varepsilon-k\varepsilon-\alpha+(\beta+k)\left(a\cot(\pi k)\right)-(\alpha\beta-\alpha)x\right),
\end{array}
\end{equation}
is smooth for all values of $a$, including $a=0$ (staying away from $k=0,k=1$).  Therefore a more standard stability analysis of the relevant point (\ref{eq:relpoint}) is possible.

To use a Taylor series expansion of the system to analyze the nearby behavior, let 
\[
\xi=x-\frac{3}{4},
\]
and 
\[
\psi=k-\frac{1}{3}.
\]

The expansion of $f_1=a(1-x-kx)$ around the point is straightforward, and yields 
\[
f_1(\xi,\psi,a,\varepsilon)=-\frac{4}{3}a\xi-\frac{3}{4}a\psi.
\]
The expansion of the $\dot{k}$ equation is more complicated because of the number of terms, but calculation shows that the Taylor expansion is 
\[
f_2(\xi,\psi,a,\varepsilon)=-\dfrac{5}{8\pi}\varepsilon+\dfrac{5}{8\pi\sqrt{3}}a+\dfrac{3}{10\pi}\xi+\left(\dfrac{\sqrt{3}}{4\pi}-\dfrac{5}{12}\right)a\psi-\left(\dfrac{5}{4\sqrt{3}}+\dfrac{3}{4\pi}\right)\varepsilon\psi+\dfrac{\sqrt{3}}{5}\xi\psi+\mathcal{O}(3).
\]

\begin{figure}[t]
	\centering
	\includegraphics[height=3in]{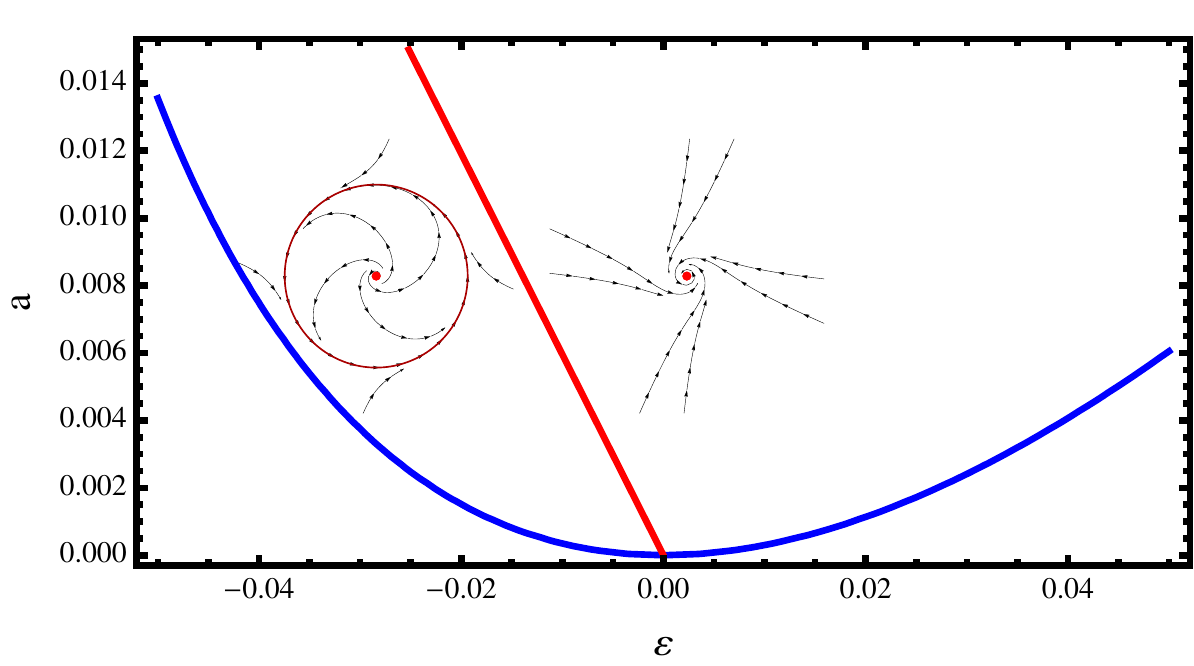}
	%\subfloat[]{\includegraphics[scale=1]{figures/ebmNonsmooth.eps}}
	%\subfloat[]{\includegraphics[width=0.45\textwidth]{figures/ebmFig.eps}}
	%\subfloat[]{\includegraphics[width=0.45\textwidth]{figures/ebmNonsmooth.eps}}
	\caption{A plot of the trace (red) and discriminant (blue) of the Jacobian matrix.  The smooth Hopf bifurcation occurs along the red curve, and limits to the nonsmooth bifurcation point.}
	\label{fig:taylor}
\end{figure}

The Jacobian of the system of equations, using the expansion, is
\begin{equation}
\label{eq:Jacobian}
J=\left[\begin{array}{cc} -\frac{4}{3}a & -\frac{3}{4} a \\ c_1 & c_3\varepsilon+c_4a\end{array}\right]=\left[\begin{array}{cc} -\frac{4}{3}a & -\frac{3}{4} a \\ \frac{3}{10\pi} & -\left(\frac{5}{4\sqrt{3}}+\frac{3}{4\pi}\right)\varepsilon+\left(-\frac{5}{12}+\frac{\sqrt{3}}{4\pi}\right)a\end{array}\right].
\end{equation}
This is singular when $a=\varepsilon=0$.  A calculation of the trace and discriminant of this Jacobian indicates the existence of a Hopf bifurcation in a neighboorhood of the singular point $a=\varepsilon=0$.  A tedious but necessary calculation shows that the genericity condition holds, and the bifurcation is supercritical for $a$ small, positive.  The trace of the jacobian is
\[
\text{tr}=\left(\frac{-21}{12}+\frac{\sqrt{3}}{4\pi}\right)a-\left(\frac{5}{4\sqrt{3}}+\frac{3}{4\pi}\right)\varepsilon,
\]
and so the trace is zero along a line in $(\varepsilon,a)$ space.  To the right of this line, the equilibrium point is stable, and to the left the equilibrium point becomes unstable, and a small stable periodic orbit forms from a supercritical Hopf bifurcation.  As $a\rightarrow0$, the bifurcation line limits to the point $(\varepsilon,a)=(0,0)$ in parameter space, showing that the supercritical Hopf bifurcation point limits to the nonsmooth fused focus bifurcation point as $a\rightarrow0$ (see Figure \ref{fig:taylor}).  It should be reiterated that the coordinate change given in equation \eqref{eq:Welblowup} is critical to this analysis, as the fused focus bifurcation here occurs in a smooth region of phase space.  This allows the Taylor expansion and pointwise limit of the supercritical Hopf bifurcation.

An astute reader may notice that the Jacobian in \eqref{eq:Jacobian} has a double zero eigenvalue at the nonsmooth bifurcation point.  This might indicate that we should expect a Takens-Bogdanov bifurcation structure, instead of a simple supercritical Hopf bifurcation.  However, this does not occur here because the parameter region $a<0$ is not allowed by the system, or, if one chooses to include it, one forces an extra symmetry as a negative sign in the inverse tangent in \eqref{eq:ksmooth}.  In the smooth system, when $a>0$, the periodic orbit which is formed by the Hopf bifurcation is later destroyed by a periodic orbit from a subcritical Hopf bifurcation. This limits to a different kind of nonsmooth bifurcation as $a\rightarrow0$, which destroys the nonsmooth periodic orbit, but which is not immediately relevant to this result, so will not be discussed in detail here.

\vspace{1em}

So, the fused focus bifurcation in Welander's nonsmooth model is a direct extension of a smooth supercritical Hopf bifurcation.

\section{Discussion}
\label{sec:discussion}

The relationship between smooth and nonsmooth systems is not mathematically obvious.  In many applications, including Welander's ocean convection model, the choice to use a nonsmooth model is made based on ease of computation, and not for physically relevant reasons.  However, it has been established that nonsmooth systems are not guaranteed to behave in a qualitatively similar way to their smooth counterparts \cite{Jeffrey11, Jeffrey14hidden}.  This work is another step towards a better understanding of this relationship.  In this model, there is an example of a nonsmooth analogue to a standard supercritical Hopf bifurcation, and it is rigorously shown that this is an extension of a standard smooth supercritical Hopf bifurcation.  This result is not dependent on choosing a way to smooth out the nonsmooth system, but is rather given by the original model itself.  Without physical reasons for the behavior of the system on the splitting manifold, one cannot conclude that the model is physically relevant there.  However, it is certain that the nonsmooth model in this case is not losing interesting behavior which might appear in the ``true" smooth model.

Mathematically, it is clear that one must use caution when using nonsmooth models to approximate smooth phenomena.  The method introduced here cannot be extended to all nonsmooth bifurcations, but it is possible that it can be used for nonsmooth bifurcations that are direct analogues to smooth phenomena, such as the fused focus bifurcation.  This illuminating example may provide inspiration for the development of a similar, more general theorem about perturbation of these bifurcation structures.

\section*{Acknowledgments}
The author would like to thank Richard McGehee, Mike Jeffrey, Paul Glendinning, and Mary Silber for their input and advice on this work.

\nocite{*}
\newpage
\bibliographystyle{amsplain}
\bibliography{researchbib}

\end{document}